\documentstyle[12pt]{amsart}
\theoremstyle{plain}
\newtheorem{Thm}{}[section]

\title {Existence of vector bundles of rank two with fixed determinant and sections}
\author{ Montserrat Teixidor i Bigas}
\address {Mathematics Department, Tufts University, Medford MA 02155}
\email {montserrat.teixidoribigas@@tufts.edu}
\begin{document}
\maketitle

\begin{abstract}
Consider the scheme $B_{2,L}^k$ of stable vector bundles of rank two
and fixed determinant $L$ which have at least $k$ sections. Under
suitable numerical conditions and for generic $L$, we show the
existence of a component of the expected dimension of $B_{2,L}^k$.

\end{abstract}

\begin{section}{Introduction}

The determinant of a vector bundle is the line bundle obtained as
its highest wedge power. Given a projective non-singular curve of
genus $g$ and a line bundle $L$ on it, there exists a moduli space
of dimension $(r^2-1)(g-1)$ parameterizing vector bundles of fixed
rank and with the given line bundle as determinant. We are
interested in the Brill-Noether locus $B^k_{r,L}$ consisting of
stable vector bundles of rank $r$ with determinant $L$ with at least
$k$ independent sections.  As in the case that one fixes only the
degree and not the determinant(see \cite{ivona} for an overview of
results in that case), this locus can be represented as a
determinatal variety. Therefore, its expected dimension  is given by
the Brill-Noether number
$$\rho ^k_{r,L}=(r^2-1)(g-1)-k(k-d+r(g-1)).$$

While there has been substantial progress in our knowledge of the
non-fixed determinant case, no work has been done so far in trying
to see how far reality is from this expectation. The purpose of this
paper is to partially fill this gap by showing existence of a
component of the expected dimension for stable vector bundles of
rank two with fixed generic determinant and a preassigned number of
sections if the degree is sufficiently large.

 For small degree $d$
or for special line bundles, things are expected to behave
differently. It is known that when the determinant is the canonical,
the locus has dimension larger than the number above (see
\cite{M1},\cite{M2},  \cite{BF}, \cite{r2detK}), \cite{Petri}.
Osserman recently extended  this set up to other special line
bundles (\cite{O}).

 \bigskip

In this paper, we deal with the case in which $r=2$ The result that
we obtain is the following

\begin{Thm} {\bf Theorem}
Let $C$ be a generic curve of genus $g$ and $L$ a generic line
bundle on $C$,  $B_{2,L}^k$ is non-empty  and has a component of the
expected dimension $\rho $ if $\rho\ge g-3$ for $k$ even. For
$k=2k_1+1$ odd, $B_{2,L}^k$ is non-empty  and has a component of
dimension the expected dimension if $(k_1+1)(k-d+2(g-1))\le g-1$.
\end{Thm}

\end{section}

\begin{section}{Review of some technical results}

The main point of the proof is the following fact (a similar
argument  was already used in \cite{duke}, \cite {Petri}): the
dimension of $B^k_{2,L}$ at any point is at least $\rho$. One can
also consider the case of a family of curves
$${\cal C}\rightarrow T$$
and a line bundle ${\cal L}$ on ${\cal C}$ of degree $d$  on every
fiber. Define
$${\cal B}^k_{2,{\cal L}}=\{ (b,E)|b\in B,\  E \in B^k_{2,{\cal L}_b}(C_b)\}.$$
Then, $$dim{\cal B}^k_{2,{\cal L}}\ge \rho ^k_{2,L}+dim B$$ at every
point. If one can find a particular curve $C_0$ such that the
dimension of $B^k_{2,L}(C_0)$ is $\rho$, then the dimension of the
generic fiber of the map ${\cal B}^k_{2,L}\rightarrow T$ attains its
minimum in a neighborhood of the curve. Hence, for a generic curve
$C$ in an open neighborhood of $C_0$, the dimension of
$B^k_{2,L}(C)$ is $\rho$ (and the locus is non-empty).
 We only need to explicitly
 exhibit a curve $C_0$ and the
corresponding family of vector bundles in $B^k_{2,L}(C_0)$. Our
$C_0$ is a reducible curve that we define as follows

\begin{Thm} {\bf Definition} \label{cce}
Let $C_1...C_g$ be elliptic curves. Let $P_i,Q_i$ be generic
points in $C_i$. Then $C_0$ is the  chain obtained by gluing the
elliptic curves when identifying the point $Q_i$ in $C_i$
 to the point $P_{i+1}$ in $C_{i+1},\ i=1...g-1$.
\end{Thm}

When dealing with reducible curves, the notion of a line bundle
and a space of its sections needs to be replaced by the analogous
concept of limit linear series as introduced by Eisenbud and
Harris (cf\cite{EH}).
 A similar definition can be given
for vector bundles (cf\cite{duke}, \cite{Clay} ). For the
convenience of the reader, we reproduce this definition here

\begin{Thm}\label{lls}
 {\bf Limit linear series}
A limit linear series of rank $r$, degree $d$ and dimension $k$ on a
chain of $M$ (not necessarily elliptic) curves consists of data I,II
below for which data III, IV exist satisfying conditions a)-c)

(I) For every component $C_i$, a vector bundle $E_i$ of rank $r$ and
degree $d_i$ and a $k$-dimensional space $V_i$ of sections of $E_i$

(II) For every node obtained by gluing $Q_i$ and $P_{i+1}$, an
isomorphism of the projectivisation of the fibers $(E_i)_{Q_i}$ and
$(E_{i+1})_{P_{i+1}}$

(III) A positive integer $b$

(IV) For every node obtained by gluing $Q_i$ and $P_{i+1}$ basis
$s^t_{Q_i}, s^t_{P_{i+1}},\ t=1...k$ of the vector spaces $V_i,
V_{i+1}$ in (I).

\medskip
Subject to the conditions

(a) $\sum_{i=1}^M d_i-r(M-1)b=d$

(b) The orders of vanishing at $Q_i, P_{i+1}$ of the sections of the
chosen basis satisfy
$ord_{Q_i}s^t_{Q_i}+ord_{P_{i+1}}s^t_{P_{i+1}}\ge b$

(c) Sections of the vector bundles $E_i(-bP_i), E_i(-bQ_i)$ are
completely determined by their value at the nodes.
\end{Thm}

{\bf Notation}

 We shall denote by
$u^i_j,\ v^i_j\ j=1...k,\ i=1...g$ the value of these vanishing. We
shall always assume that the $u^i_j$ are in increasing order and the
$v^i_j$ are in decreasing order, namely
$$u^i_1\le u^i_2\le ...\le u^i_k,\ v^i_1\ge v^i_2\ge ...\ge v^i_k$$

At most two of the $u^i$ can take a given value (say
$u^i_{j-1}<u^i_j=u^i_{j+1}<u^i_{j+2}$ and when this happens, there
are two linearly independent sections of $E_i(-u^i_jP_i)$ that
vanish at $P_i$ with multiplicity exactly $u^i_j$ and generate the
fiber of $E_i(-u^i_jP_i)$  at this point. Analogous statements can
be given for the vanishing at $Q_i$.

\begin {Thm} \label{notacio} {\bf Notation.}\end{Thm}
We shall write
$$d=2d_1+\epsilon, \ k=2k_1+\delta,\ 0\le \epsilon ,\ \delta \le 1$$
We shall use $b=d_1$.

To simplify notations, we shall omit the superindex $i$ when the
curve is clear. In order to optimize the vanishing, we always take
the $u^1_j$ at $P_1$ to be the smallest possible, namely
$(0,0,1,1,2,2...)$ and for $i>1$, take
$u^i_j=b-v^{i-1}_j=d_1-v^{i-1}_j$ where $b$ is as defined in
\ref{lls} III.

\begin {Thm} \label{estable} {\bf Remark.}\end{Thm}
We need to prove that some vector bundles we construct on the
reducible curves are stable. From \cite{arbre} Step 2, p.342 and
\ref{estable} Prop 1.2 it is enough to see that the restrictions to
each of the components is semistable and at least one of them is
stable or in the case they are all  strictly semistable, the
destabilizing subbundles do not glue with each other.

\bigskip
\end{section}

\begin{section}{Even degree and number of sections}
We start by considering the case of even degree and even number of
sections. So
$$d=2d_1, \ k=2k_1$$
We define $$\alpha =k_1-d_1+g-1.$$

The condition $\rho \ge g-3$ can then be written as
$2k_1(k_1-d+g-1)\le g$ or equivalently
$$(*) \ \ 2k_1\alpha \le g$$

We construct a vector bundle on a curve $C_0$ as in \ref{cce} with
fixed determinant and a $k$-dimensional limit linear series.

Define, for $i>1$
$$u^i_j=d_1-v^{i-1}_j, j=1...k$$

On the  curve $C_1$, take the vector bundle
$$({\mathcal O}(d_1Q_1))\oplus L_1$$ where $L_1$ has been chosen so that the determinant is as preassigned
(and by our genericity assumption for the determinant, this implies
it is generic).

On the  curve $C_i, i=k_1t+j,\ j=2...k_1, t=0,\cdots 2\alpha-2$,
 take the vector bundle to be
$$({\cal O}(u^i_{2j-1}P_i+(d_1-u^i_{2j-1})Q_i))\oplus L_i$$
 where $L_i$ has been chosen so that the determinant is as preassigned.
 Glue this curve to the previous one so that $L_i$ glues with the
 direction of the section on $C_{i-1}$ with vanishing
 $v^{i-1}_{2j-2}$ at $Q_{i-1}$.
Then
 $$v^i_{t}=d_1-u^i_{t}-1,\ t\not= 2j-1$$
and
$$v^i_{2j-1}=d_1-u^i_{2j-1}.$$

On the curve  $C_i, i={k_1t+1}, t=1,\cdots 2\alpha-2$,
 take the vector bundle to be the unique indecomposable vector bundle of rank two of degree
 $2d_1+1$ with preassigned  determinant if $t$ is odd and of degree
  $2d_1-1$ with preassigned determinant if $t$ is even (see \cite{A} for results on vector bundles on elliptic curves).
  This vector bundle has a unique section that vanishes at $P_i$ with
  multiplicity $u^i_1$ and at $Q_i$ with multiplicity $d_1-\epsilon
  -u^i_1$ where $\epsilon =1$ if $t$ is even, $\epsilon =0$ if $t$ is odd. It has a second
  section that vanishes at $P_i$ with
  multiplicity $u^i_{2k_1}$ and at $Q_i$ with multiplicity $d_1-\epsilon
  -u^i_{2k_1}$ where $\epsilon =1$ if $t$ is even. Glue these two
  sections with the directions of the sections on $C_{i-1}$ that
  vanish to order $v^{i-1}_1$ and $v^{i-1}_{2k_1}$ respectively.
  Then if $t$ is even or $t=1$
   $$v^i_{t}=d_1-u^i_{t}-\epsilon$$
and for odd $t\not= 1$
$$v^i_{t}=d_1-u^i_t-1.$$

  On the  curve $C_i, i={k_1(2\alpha-1)+j},\ j=1...k_1-1$,
 take the vector bundle to be
$$({\cal O}(u^i_{2j}P_i+(d_1-u^i_{2j})Q_i))\oplus L_i$$
 where $L_i$ has been chosen so that the determinant is as preassigned.
 Glue this curve to the previous one so that $L_i$ glues with the
 direction of the section on $C_{i-1}$ with vanishing
 $v^{i-1}_{2j-1}$ at $Q_{i-1}$.
Then
 $$v^i_{t}=d_1-u^i_{t}-1,\ t\not= 2j$$
and
$$v^i_{2j}=d_1-u^i_{2j}.$$

 On the  curve $C_i, i={k_1(2\alpha-1)+k_1=2k_1\alpha}$,
 take the vector bundle to be
$$({\cal O}(u^i_{2k_1}P_i+(d_1-u^i_{2k_1})Q_i))\oplus L_i$$
 where $L_i$ has been chosen so that the determinant is as preassigned.
 Glue this curve to the previous one so that $L_i$ glues with the
 direction of the section on $C_{i-1}$ with vanishing
 $v^{i-1}_{2k_1-1}$ at $Q_{i-1}$ and $({\cal O}(u^i_{2k_1}P_i+(d_1-u^i_{2k_1})Q_i))$
  glues with the direction with vanishing
 $v^{i-1}_{2k_1}$ at $Q_{i-1}$.
Then
 $$v^i_{t}=d_1-u^i_{t}-1,\ t\not= 2k_1$$
and
$$v^i_{2k_1}=d_1-u^i_{2k_1}.$$

Note that, by our assumptions (se (*)), $g-2k_1\alpha \ge 0$.

On the remaining $g-2k_1\alpha $ components, take the vector bundle
to be the direct sum of two line bundles of degree $d_1$ such that
their tensor product is the preassigned determinant. Take the gluing
among these components to be generic
$$v^i_j=d_1-u^i_j-1$$
Hence,
$$(v^g_1,v^g_2,....v^g_{k-1},v^g_k)=(d_1-g+\alpha, d_1-g+\alpha,...
d_1-g+\alpha-(k_1-1),d_1-g+\alpha-(k_1-1))$$
$$=(k_1-1,k_1-1,...0,0)$$
The vanishings at $Q_g$ are the smallest possible. Hence,
 we cannot make the vector bundles
more general and still obtain a linear series of dimension $k$.

We now compute the number of moduli of such a family. The
restrictions of the vector bundles to the first $2k_1\alpha$
components are completely determined. On the remaining
$g-2k_1\alpha$ components, the restriction depends on one parameter.
The  gluing at each of the last $g-2k_1\alpha$ nodes is generic and
therefore depends on four parameters. At the nodes $k_1t,\ 1\le t\le
2\alpha-2$, they depend on two parameters while at the remaining
nodes, they depend on three parameters. Each of the vector bundles
obtained as the restriction to a component has a two dimensional
family of automorphisms except for  the vector bundles on the
components $k_1t+1,\ 1\le t\le 2\alpha-2$  which have only one. The
resulting vector bundle on the reducible curve has a one-dimensional
family of automorphisms, as it is stable. Hence, the number of
moduli for the family is
$$g-2k_1\alpha +4(g-2k_1\alpha)+2(2\alpha-2+1)+3( 2k_1\alpha-(2\alpha-2)-2)-2(g
-(2\alpha-2))-(2\alpha-2)+1=\rho$$

\end{section}

\begin{section}{Odd degree and even number of sections}
 Now
$$d=2d_1+1, \ k=2k_1$$
We define $$ \alpha =k_1-d_1+g-1.$$

Now the condition $\rho \ge g-3$ can then be written as
$$(*) \ \ k_1(2\alpha -1)\le g.$$
 This case is very similar to the
previous one, except that we include an odd number of vector bundles
of odd degree.

On the  curve $C_1$, take the vector bundle
$$({\mathcal O}(d_1Q_1))\oplus L_1$$ where $L_1$ has been chosen so that the determinant is as preassigned
(and by our genericity assumption for the determinant, this implies
it is generic).

On the  curve $C_i, i=k_1t+j,\ j=2...k_1, t=0,\cdots 2\alpha-3$,
 take the vector bundle to be
$$({\cal O}(u^i_{2j-1}P_i+(d_1-u^i_{2j-1})Q_i))\oplus L_i$$
 where $L_i$ has been chosen so that the determinant is as preassigned.
 Glue this curve to the previous one so that $L_i$ glues with the
 direction of the section on $C_{i-1}$ with vanishing
 $v^{i-1}_{2j-2}$ at $Q_{i-1}$.
Then
 $$v^i_{l}=d_1-u^i_{l}-1,\ l\not= 2j-1$$
and
$$v^i_{2j-1}=d_1-u^i_{2j-1}.$$

On the curve  $C_i, i={k_1t+1}, t=1,\cdots 2\alpha-1$,
 take the vector bundle to be the unique indecomposable vector bundle of rank two and degree
 $2d_1+1$ with preassigned  determinant if $t$ is odd and of degree
  $2d_1-1$ with preassigned determinant if $t$ is even.
  This vector bundle has a unique section that vanishes at $P_i$ with
  multiplicity $u^i_1$ and at $Q_i$ with multiplicity $d_1-\epsilon
  -u^i_1$ where $\epsilon =1$ if $t$ is even, $\epsilon =0$ if $t$ is odd. It has a second
  section that vanishes at $P_i$ with
  multiplicity $u^i_{2k_1}$ and at $Q_i$ with multiplicity $d_1-\epsilon
  -u^i_{2k_1}$ where $\epsilon =1$ if $t$ is even. Glue these two
  sections with the directions of the sections on $C_{i-1}$ that
  vanish to order $v^{i-1}_1$ and $v^{i-1}_{2k_1}$ respectively.
  Then if $l$ is even or $l=1$
   $$v^i_{l}=d_1-u^i_{l}-\epsilon$$
and for odd $l\not= 1$
$$v^i_{l}=d_1-u^i_l-1.$$

  On the  curve $C_i, i={k_1(2\alpha-2)+j},\ j=1...k_1-1$,
 take the vector bundle to be
$$({\cal O}(u^i_{2j}P_i+(d_1-u^i_{2j})Q_i))\oplus L_i$$
 where $L_i$ has been chosen so that the determinant is as preassigned.
 Glue this curve to the previous one so that $L_i$ glues with the
 direction of the section on $C_{i-1}$ with vanishing
 $v^{i-1}_{2j-1}$ at $Q_{i-1}$.
Then
 $$v^i_{t}=d_1-u^i_{t}-1,\ t\not= 2j$$
and
$$v^i_{2j}=d_1-u^i_{2j}.$$

 On the  curve $C_i, i=k_1(2\alpha-2)+k_1=k_1(2\alpha -1)$,
 take the vector bundle to be
$$({\cal O}(u^i_{2k_1}P_i+(d_1-u^i_{2k_1})Q_i))\oplus L_i$$
 where $L_i$ has been chosen so that the determinant is as preassigned.
 Glue this curve to the previous one so that $L_i$ glues with the
 direction of the section on $C_{i-1}$ with vanishing
 $v^{i-1}_{2k_1-1}$ at $Q_{i-1}$ and $({\cal O}(u^i_{2j}P_i+(d_1-u^i_{2j})Q_i))$
  glues with the direction with vanishing
 $v^{i-1}_{2k_1}$ at $Q_{i-1}$.
Then
 $$v^i_{l}=d_1-u^i_{l}-1,\ t\not= 2k_1$$
and
$$v^i_{2k_1}=d_1-u^i_{2k_1}.$$

Note that by our assumptions (see (*)), $g- k_1(2\alpha -1)\ge 0$.
On the remaining $g-k_1(2\alpha -1)$ components, take the vector
bundle to be the direct sum of two line bundles of degree $d_1$ such
that their tensor product is the preassigned determinant. Take the
gluing among these components to be generic
$$v^i_j=d_1-u^i_j-1$$
Hence,
$$(v^g_1,v^g_2,....v^g_{k-1},v^g_k)=(d_1-g+\alpha, d_1-g+\alpha,...
d_1-g+\alpha-(k_1-1),d_1-g+\alpha-(k_1-1))$$
$$=(k_1-1,k_1-1,...0,0)$$
The vanishings at $Q_g$ are the smallest possible. Hence,
 we cannot make the vector bundles
more general and still obtain a linear series of dimension $k$.

We now compute the number of moduli of such a family. The
restrictions of the vector bundles to the first $k_1(2\alpha-1)$
components are completely determined. On the remaining
$g-k_1(2\alpha -1)$ components, the restriction depends on one
parameter. The  gluing at each of the last $g-k_1(2\alpha-1)$ nodes
is generic and therefore depends on four parameters. At the nodes
$k_1t,\ 1\le t\le 2\alpha-3$ or the node $k_1(2\alpha -1)-1$, they
depend on two parameters while at the remaining nodes, they depend
on three parameters. Each of the vector bundles obtained as the
restriction to a component has a two dimensional family of
automorphisms except for  the vector bundles on the components
$k_1t+1,\ 1\le t\le 2\alpha-3$  which have only one. The resulting
vector bundle on the reducible curve has a one-dimensional family of
automorphisms, as it is stable. Hence, the number of moduli for the
family is
$$g-k_1(2\alpha-1) +4(g-k_1(2\alpha) -1)+2(2\alpha-2)+3( k_1(2\alpha-1)-$$
$$-(2\alpha-2)-1)-2(g-2(g-(2\alpha-3))-(2\alpha-3)+1=\rho$$

\end{section}

\begin{section}{Even degree and odd number of sections}

Write
$$d=2d_1, \ k=2k_1+1,\ \alpha =k_1-d_1+g-1.$$

On the curve  $C_i, i={(k_1+1)t+1}, t=0,\cdots 2\alpha +1$,
 take the vector bundle to be the unique indecomposable vector bundle of rank two and degree
 $2d_1+1$ with preassigned  determinant if $t$ is even and of degree
  $2d_1-1$ with preassigned determinant if $t$ is odd.
  This vector bundle has a unique section that vanishes at $P_i$ with
  multiplicity $u^i_{2k_1+1}$ and at $Q_i$ with multiplicity $d_1-\epsilon
  -u^i_{2k_1}$ where $\epsilon =1$ if $t$ is odd. Glue this
  section with the directions of the sections on $C_{i-1}$ that
  vanishes to order  $v^{i-1}_{2k_1+1}$.
  Then if $l$ is odd
   $$v^i_{l}=d_1-u^i_{l}-\epsilon$$
and for even $l$
$$v^i_{l}=d_1-u^i_l-1-\epsilon.$$

On the  curve $C_i, i=(k_1+1)t+j,\ j=2...k_1+1, t=0,\cdots 2\alpha$,
 take the vector bundle to be
$$({\cal O}(u^i_{2j-2}P_i+(d_1-u^i_{2j-2})Q_i))\oplus L_i$$
 where $L_i$ has been chosen so that the determinant is as preassigned.
 Glue this curve to the previous one so that $L_i$ glues with the
 direction of the section on $C_{i-1}$ with vanishing
 $v^{i-1}_{2j-3}$ at $Q_{i-1}$.
Then
 $$v^i_{l}=d_1-u^i_{l}-1,\ l\not= 2j-2$$
and
$$v^i_{2j-2}=d_1-u^i_{2j-1}.$$

The condition that we are assuming, namely $(k_1+1)(k-d+2(g-1))\le
g-1$ can be written as $g-(k_1+1)(2\alpha+1)-1\ge 0$.

On the remaining $g-(k_1+1)(2\alpha +1)-1$ components, take the
vector bundle to be the direct sum of two  line bundles of degree
$d_1$ such that their tensor product is the preassigned determinant.
Take the gluing among these components to be generic. Then,
$$v^i_j=d_1-u^i_j-1$$
Hence,
$$(v^g_1,v^g_2,....v^g_{k-1},v^g_k)=(d_1-g+\alpha, d_1-g+\alpha,...
d_1-g+\alpha-(k_1-1),d_1-g+\alpha-(k_1-1))$$
$$=(k_1-1,k_1-1,...0,0)$$
The vanishings at $Q_g$ are the smallest possible. Hence,
 we cannot make the vector bundles or gluing
more general and still obtain a linear series of dimension $k$.

We now compute the number of moduli of such a family. The
restrictions of the vector bundles to the first
$(k_1+1)(2\alpha+1)+1$ components are completely determined. On the
remaining $g-(k_1+1)(2\alpha +1) -1)$ components, the restriction
depends on one parameter. The  gluing at each of the last
$g-(k_1+1)(2\alpha +1)-1)$ nodes is generic and therefore depends on
four parameters. At the
 remaining nodes  it
depends on three parameters. Each of the vector bundles obtained as
the restriction to a component has a two dimensional family of
automorphisms except for  the vector bundles on the components
$(k_1+1)t+1,\ 0\le t\le 2\alpha+1$  which have only one. The
resulting vector bundle on the reducible curve has a one-dimensional
family of automorphisms, as it is stable (see \ref{estable}. Hence,
the number of parameters for the family is
$$g-(k_1+1)(2\alpha+1)-1 +4(g-(k_1+1)(2\alpha +1)) -1)+3((k_1+1)(2\alpha+1)-$$
$$-2(g-2\alpha-2)-(2\alpha+2)+1=\rho$$

\end{section}

\begin{section}{Odd degree and odd number of sections}

This case is similar to the previous one. Write
$$d=2d_1+1, \ k=2k_1+1,\ \alpha =k_1-d_1+g-1.$$

On the curve  $C_i, i={(k_1+1)t+1}, t=0,\cdots 2\alpha$,
 take the vector bundle to be the unique indecomposable vector bundle of rank two and degree
 $2d_1+1$ with preassigned  determinant if $t$ is even and of degree
  $2d_1-1$ with preassigned determinant if $t$ is odd.
  This vector bundle has a unique section that vanishes at $P_i$ with
  multiplicity $u^i_{2k_1+1}$ and at $Q_i$ with multiplicity $d_1-\epsilon
  -u^i_{2k_1}$ where $\epsilon =1$ if $t$ is odd. Glue this
  section with the directions of the sections on $C_{i-1}$ that
  vanishes to order  $v^{i-1}_{2k_1+1}$.
  Then if $l$ is odd
   $$v^i_{l}=d_1-u^i_{l}-\epsilon$$
and for even $l$
$$v^i_{l}=d_1-u^i_l-1-\epsilon.$$

On the  curve $C_i, i=(k_1+1)t+j,\ j=2...k_1+1, t=0,\cdots
2\alpha-1$,
 take the vector bundle to be
$$({\cal O}(u^i_{2j-2}P_i+(d_1-u^i_{2j-2})Q_i))\oplus L_i$$
 where $L_i$ has been chosen so that the determinant is as preassigned.
 Glue this curve to the previous one so that $L_i$ glues with the
 direction of the section on $C_{i-1}$ with vanishing
 $v^{i-1}_{2j-3}$ at $Q_{i-1}$.
Then
 $$v^i_{l}=d_1-u^i_{l}-1,\ l\not= 2j-2$$
and
$$v^i_{2j-2}=d_1-u^i_{2j-1}.$$

The condition that we are assuming, namely $(k_1+1)(k-d+2(g-1))\le
g-1$ can be written as $g-(k_1+1)2\alpha-1\ge 0$.

On the remaining $g-(k_1+1)2\alpha -1$ components, take the vector
bundle to be the direct sum of two  line bundles of degree $d_1$
such that their tensor product is the preassigned determinant. Take
the gluing among these components to be generic. Then,
$$v^i_j=d_1-u^i_j-1$$
Hence,
$$(v^g_1,v^g_2,....v^g_{k-1},v^g_k)=(d_1-g+\alpha, d_1-g+\alpha,...
d_1-g+\alpha-(k_1-1),d_1-g+\alpha-(k_1-1))$$
$$=(k_1-1,k_1-1,...0,0)$$
The vanishings at $Q_g$ are the smallest possible. Hence,
 we cannot make the vector bundles or gluing
more general and still obtain a linear series of dimension $k$.

We now compute the number of moduli of such a family. The
restrictions of the vector bundles to the first $(k_1+1)2\alpha+1$
components are completely determined. On the remaining
$g-(k_1+1)2\alpha -1$ components, the restriction depends on one
parameter. The  gluing at each of the last $g-(k_1+1)2\alpha -1$
nodes is generic and therefore depends on four parameters. At the
 remaining nodes  it
depends on three parameters. Each of the vector bundles obtained as
the restriction to a component has a two dimensional family of
automorphisms except for  the vector bundles on the components
$(k_1+1)t+1,\ 0\le t\le 2\alpha$  which have only one. The resulting
vector bundle on the reducible curve has a one-dimensional family of
automorphisms, as it is stable (see \ref{estable}). Hence, the
number of parameters for the family is
$$g-(k_1+1)2\alpha-1 +4(g-(k_1+1)(2\alpha) -1)+3((k_1+1)(2\alpha)-$$
$$-2(g-(2\alpha+1))-(2\alpha+1)+1=\rho$$

\end{section}


\bigskip
\begin{thebibliography}{cccc}
\bibitem [A]{A} M.Atiyah, {\it Vector Bundles over an elliptic
curve}, Proc.London Math.Soc. {\bf 3} (1957), 414-452.

\bibitem[BF]{BF} A.Bertram, B.Feinberg {On stable rank two bundles
with canonical determinant and many sections}, in "Algebraic
Geometry", ed P.Newstead, Marcel Dekker 1998, 259-269.

\bibitem [EH]{EH} D.Eisenbud, J.Harris {\it Limit linear series,
basic theory}, Invent.Math. {\bf 85} (1986), 337-371.

\bibitem [GT]{ivona} I.Grzegorczyk, M.Teixidor {\it Brill Noether Theory for
stable vector bundles} in "Moduli spaces of Vector Bundles"  London
Mathematical Society Lecture Note Series {\bf 359} Cambridge
University Press 2009, 29-50.

\bibitem[M1]{M1} S.Mukai, {\it Curves and Brill-Noether Theory},
MSRI publications {\bf 28}, Cambridge University Press 1995,
145-158.

\bibitem[M2]{M2} S.Mukai {\it Non-abelian Brill-Noether Theory and
Fano threefolds}.  Sugaku Expositions  {\bf 14}  (2001),  no. 2,
125--153.

\bibitem[O]{O} B.Osserman, {\it Brill-Noether loci for fixed
determinant in rank two} arxiv 1005.0448.

\bibitem [T1]{duke} M.Teixidor {\it Brill Noether Theory for
stable vector bundles}, Duke Math. J. {\bf 62 \ N2}, (1991) 385-400

\bibitem [T2]{arbre} M.Teixidor {\it Moduli spaces of semistable
vector bundles on tree-like curves}, Math.Ann. {\bf 290} (1991).
341-348.

\bibitem [T3]{reduible} M.Teixidor {\it Moduli spaces of vector bundles on
reducible curves}, Amer J of Math. {\bf 219}(1995), 477-484.

\bibitem [T4]{r2detK} M.Teixidor {\it Rank two vector bundles with
canonical determinant.} Math. Nachr. {\bf 265} (2004), 100-106.

\bibitem [T5]{Petri} M.Teixidor {\it Petri map for rank two vector bundles with
canonical determinant.} Comp. Math. {\bf } (200), 10-10.

\bibitem [T6]{Clay} M.Teixidor {\it Vector Bundles on reducible curves and applications}
Clay mathematics Institute Proceedings, to appear.



\end{thebibliography}
\end{document}